\documentclass{amsart}

\usepackage{amsmath} 
\usepackage{amssymb}

\newtheorem{theorem}{Theorem}[section] 

\newtheorem{lemma}[theorem]{Lemma} 
\newtheorem{observation}[theorem]{Observation}

\theoremstyle{definition}
\newtheorem{definition}[theorem]{Definition}

\theoremstyle{remark}

\newcommand{\conc}{{}^\frown\!}
\newcommand{\lh}{{\ell g}} 
\newcommand{\lbd}{{\leq_{J^{\rm bd}_\omega}}}
\newcommand{\vare}{\varepsilon}
\newcommand{\vtl}{\vartriangleleft}

\newcommand{\cA}{{\mathcal A}}

\newcommand{\gb}{{\mathfrak b}}

\newcommand{\cF}{{\mathcal F}}

\newcommand{\gf}{{\mathfrak f}}

\newcommand{\frakg}{{\mathfrak g}}

\newcommand{\gh}{{\mathfrak h}}
\newcommand{\cI}{{\mathcal I}}

\newcommand{\cP}{{\mathcal P}}

\newcommand{\cT}{{\mathcal T}}

\newcommand{\cU}{{\mathcal U}}

\begin{document}

\title{Groupwise density cannot be much bigger than the unbounded number}  

\author{Saharon Shelah}
\address{Institute of Mathematics\\
 The Hebrew University of Jerusalem\\
 91904 Jerusalem, Israel\\
 and  Department of Mathematics\\
 Rutgers University\\
 New Brunswick, NJ 08854, USA}
\email{shelah@math.huji.ac.il}
\urladdr{http://shelah.logic.at}

\subjclass{Primary 03E17; Secondary: 03E05, 03E20}
\date{August 2007}

\begin{abstract}
We prove that $\frakg$ (the groupwise density  number) is smaller or equal to
$\gb^+$, the successor of the minimal cardinality of an unbounded subset of
${}^\omega \omega$. 
\end{abstract}

\maketitle

\section{Introduction}
In the present note we are interested in two cardinal characteristics of the
continuum, the unbounded number $\gb$ and the groupwise density number
$\frakg$. The former cardinal belongs to the oldest and most studied
cardinal invariants of the continuum (see, e.g., van Douwen \cite{vD} and
Bartoszy\'nski and Judah \cite{BaJu95}) and it is defined as follows. 

\begin{definition}
\begin{enumerate}
\item[(a)]  The partial order $\lbd$ on ${}^\omega\omega$ is defined by  
\[f\lbd g  \mbox{ if and only if } (\exists N<\omega)(\forall n>N)(f(n)
\leq g(n)).\] 
\item[(b)]  The {\em unbounded number }$\gb$ is defined by  
\[\gb=\min\{|\cF|:\cF \subseteq {}^\omega\omega\mbox{ has no $\lbd$--upper
  bound in } {}^\omega\omega\}.\]
\end{enumerate}
\end{definition}

The groupwise density number $\frakg$, introduced in Blass and Laflamme
\cite{BsLa89}, is perhaps less popular but it has gained substantial
importance in the realm of cardinal invariants. For instance, it has been
studied in connection with the cofinality $c({\rm Sym}(\omega))$ of the
symmetric group on the set $\omega$ of all integers, see Thomas \cite{Th98a}
or Brendle and Losada \cite{BrLo03}. The cardinal $\frakg$ is defined as
follows.

\begin{definition}
\begin{enumerate}
\item[(a)] We say that a family $\cA\subseteq [\omega]^{\aleph_0}$ is {\em
    groupwise dense}  whenever:  
\begin{itemize}
\item $B \subseteq A \in \cA$, $B \in [\omega]^{\aleph_0}$ implies $B\in
  \cA$, and 
\item for every increasing sequence $\langle m_i:i<\omega\rangle
\in {}^\omega \omega$ there is an infinite set $\cU\subseteq\omega$ such
that $\bigcup\{[m_i,m_{i+1}):i \in \cU\}\in\cA$.
\end{itemize}
\item[(b)] The {\em groupwise density number\/} $\frakg$ is defined as the
  minimal cardinal $\theta$ for which there is a sequence $\langle
  \cA_\alpha:\alpha < \theta\rangle$ of groupwise dense subsets of
  $[\omega]^{\aleph_0}$ such that 
\[\big(\forall B \in [\omega]^{\aleph_0}\big)\big(\exists\alpha<\theta\big)
\big(\forall A \in \cA_\alpha\big)\big(B\not\subseteq^* A\big).\]
\end{enumerate}
\end{definition}
(Recall that for infinite sets $A$ and $B$,  $A \subseteq^* B$ means $A
\setminus B$ is finite.)

The unbounded number $\gb$ and  groupwise density number $\frakg$ can be in
either order, see Blass \cite{Bs89} and Mildenberger and Shelah
\cite{MdSh:731}, \cite{MdSh:843}. However, as we show in Theorem
\ref{thmbg}, $\frakg$ cannot be bigger than $\gb^+$.
\medskip

\noindent{\bf Notation}: Our notation is rather standard and compatible with
that of classical textbooks on Set Theory (like Bartoszy\'nski and Judah 
\cite{BaJu95}). We will keep the following rules concerning the use of
symbols.
\begin{enumerate}
\item $A,B,\cU$ (with possible sub- and superscripts) denote subsets of
  $\omega$, infinite if not said otherwise.
\item $m,n,\ell,k,i,j$ are natural numbers.
\item $\alpha,\beta,\gamma,\delta,\varepsilon,\xi,\zeta$ are ordinals,
  $\theta$ is a cardinal. 
\end{enumerate}

\section{The result}

\begin{lemma}
\label{lemm}
For some cardinal $\theta\leq\gb$ there is a sequence $\langle
B_{\zeta,t}:\zeta < \theta,\ t \in I_\zeta\rangle$ such that: 
\begin{enumerate}
\item[(a)] $B_{\zeta,t} \in [\omega]^{\aleph_0}$
\item[(b)] if $\zeta<\theta$ and $s \ne t$ are from $I_\zeta$, then
  $B_{\zeta,s} \cap B_{\zeta,t}$ is finite (so $|I_\zeta| \le
  2^{\aleph_0}$),  
\item[(c)] for every $B \in [\omega]^{\aleph_0}$ the set 
\[\big\{(\zeta,t):\zeta<\theta\ \&\ t\in I_\zeta\ \&\  B_{\zeta,t} \cap
B\mbox{ is infinite }\big\}\] 
is of cardinality $2^{\aleph_0}$.
\end{enumerate}
\end{lemma}

\begin{proof}
This is a weak version of the celebrated base-tree theorem of Bohuslav
Balcar and Petr Simon with $\theta=\gh$ which is known to be $\leq\gb$, see
Balcar and Simon \cite[3.4, pg.350]{BaSi89}. However, for the sake of
completeness of our exposition, let us present the proof.  

Let $\langle f_\zeta:\zeta < \gb\rangle$ be a $\lbd$--increasing sequence of
members of ${}^\omega \omega$ with no $\lbd$--upper bound in ${}^\omega
\omega$. Moreover we demand that each $f_\zeta$ is increasing.  By induction
on $\zeta<\gb$ choose sets $\cT_\zeta$ and systems $\langle
B_{\zeta,\eta}:\eta \in \cT_{\zeta +1}\rangle$ such that:
\begin{enumerate}
\item[(i)]  $\cT_\zeta \subseteq {}^\zeta(2^{\aleph_0})$ and if $\eta \in
\cT_{\zeta +1}$ then $B_{\zeta,\eta} \in [\omega]^{\aleph_0}$, 
\item[(ii)]  if $\eta \in \cT_\zeta$ and  $\varepsilon < \zeta$, then $\eta
\restriction \varepsilon \in \cT_\varepsilon$, 
\item[(iii)] if $\zeta$ is a limit ordinal, then 
\[\cT_\zeta=\big\{\eta \in {}^\zeta(2^{\aleph_0}):\big( \forall\vare < \zeta
\big)\big(\eta\restriction\vare\in \cT_\vare\big)\mbox{ and }\big(\exists A
\in [\omega]^{\aleph_0}\big)\big(\forall \varepsilon < \zeta\big)\big(A
\subseteq^*B_{\vare,\eta\restriction (\varepsilon +1)}\big)\},\]  
\item[(iv)] if $\vare<\zeta$ and $\eta \in \cT_{\zeta +1}$, then
$B_{\zeta,\eta} \subseteq^* B_{\vare,\eta \restriction (\vare+1)}$, 
\item[(v)]  for $\eta\in\cT_{\zeta +1}$ and $m_1 < m_2$ from
$B_{\zeta,\eta}$ we have $f_\zeta(m_1) < m_2$,
\item[(vi)]  if $\eta\in\cT_\varepsilon$, then the set
  $\{B_{\varepsilon,\nu}:\eta\vtl \nu \in \cT_{\vare+1}\}$ is an infinite  
maximal subfamily of  
\[\big\{A \in [\omega]^{\aleph_0}: \big(\forall\xi<\varepsilon\big)\big(A
\subseteq^* B_{\xi,\eta \restriction (\xi +1)}\big)\big\}\] 
consisting of pairwise almost disjoint sets.
\end{enumerate}
It should be clear that the choice is possible. Note that for some limit 
$\zeta<\gb$ we may have $\cT_\zeta=\emptyset$  (and then also $\cT_\xi= 
\emptyset$ for $\xi>\zeta$). Also, if we define $\cT_\gb$ as in (iii), then  
it will be empty (remember clause (v) and the choice of $\langle
f_\zeta:\zeta < \gb\rangle$). 

The lemma will readily follow from the following fact.
\begin{enumerate}
\item[$(\circledast)$]  For every $A \in [\omega]^{\aleph_0}$ there is
  $\xi<\gb$ such that
\[|\big\{\eta\in\cT_{\xi+1}: B_{\xi,\eta} \cap A\mbox{ is infinite }\big\}|
= 2^{\aleph_0}.\] 
\end{enumerate}
To show $(\circledast)$ let $A\in [\omega]^{\aleph_0}$ and define 
\[S=\bigcup_{\zeta<\gb}\big\{\eta \in \cT_\zeta: (\forall\vare<
\zeta)( A\cap B_{\vare,\eta \restriction (\vare+1)} \mbox{ is infinite
})\big\}.\] 
Clearly $S$ is closed under taking the initial segments and
$\langle\rangle\in S$.  By the ``maximal'' in clause (vi), we have that  
\begin{enumerate}
\item[$(\circledast)_1$] if  $\eta\in S \cap\cT_\zeta$ where $\zeta<\gb$ is
  non-limit or ${\rm cf}(\zeta)=\aleph_0$, \\
 then $(\exists\nu)(\eta \vtl \nu \in \cT_{\zeta +1} \cap S)$. 
\end{enumerate}
Now,
\begin{enumerate}
\item[$(\circledast)_2$] if $\eta \in S$ and $\lh(\eta)$ is non-limit or
  ${\rm cf}(\lh(\eta))=\aleph_0$, then there are $\vtl$--incomparable
  $\nu_0,\nu_1 \in  S$ extending $\eta$, i.e.,  $\eta\vtl\nu_0$ and
  $\eta\vtl\nu_1$.   
\end{enumerate}
[Why?  As otherwise $S_\eta = \{\nu\in S:\eta \trianglelefteq \nu\}$ is
linearly ordered by $\vtl$, so let $\rho = \bigcup S_\eta$. It follows from
$(\circledast)_1$ that $\lh(\rho)>\lh(\eta)$ is a limit ordinal (of
uncountable cofinality). Moreover, by (iv)+(vi), we have that  
\[\vare< \lh(\rho)\quad \Rightarrow\quad A \cap B_{\lh(\eta),\rho
  \restriction (\lh(\eta)+1)} =^* A \cap B_{\vare,\rho\restriction(\vare+
  1)}.\]  
Hence, by (iii)+(ii), $\rho\in \cT_{\lh(\rho)}$ so necessarily
$\lh(\rho)<\gb$.  Using (vi) again we may conclude that there is $\rho'\in
S$ properly extending $\rho$,  getting a contradiction.]

Consequently, we may find a system $\langle \eta_\rho: \rho \in {}^{\omega 
  >} 2\rangle\subseteq S$ such that for every $\rho \in {}^{\omega >}2$:
\begin{itemize}
\item $k<\lh(\rho)\quad \Rightarrow\quad \eta_{\rho \restriction k} \vtl
  \eta_\rho$, and  
\item $\eta_{\rho \conc \langle 0\rangle},\eta_{\rho \conc \langle
    1\rangle}$ are $\vtl$--incomparable.  
\end{itemize}
For $\rho \in {}^{\omega >}2$ let $\zeta(\rho)=\sup\{\lh(\eta_\nu): \rho 
\trianglelefteq \nu \in {}^{\omega >}2\}$. Pick $\rho$ such that
$\zeta(\rho)$ is the smallest possible (note that ${\rm
  cf}(\zeta(\rho))=\aleph_0$). Now it is possible to choose a perfect
subtree $T^*$ of ${}^{\omega >}2$ such that  
\[\nu \in \lim(T^*)\quad \Rightarrow\quad \sup\{\lh(\eta_{\nu 
  \restriction n}):n < \omega\}=\zeta(\rho).\]
We finish by noting that for every $\nu \in \lim(T^*)$ we have that
$\bigcup\{\eta_{\nu\restriction n}:n < \omega\}\in 
\cT_{\zeta(\rho)}\cap S$ and there is $\eta^*\in \cT_{\zeta(\rho)+1}\cap S$ 
extending $\bigcup\{\eta_{\nu\restriction n}:n < \omega\}$.   
\end{proof}

\begin{theorem}
\label{thmbg}
  $\frakg \leq \gb^+$.
\end{theorem}

\begin{proof}
Assume towards contradiction that  $\frakg>\gb^+$.

Let $\langle f_\alpha:\alpha<\gb\rangle\subseteq {}^\omega\omega$ be an
$\lbd$--increasing sequence with no $\lbd$--upper bound. We also demand that
all functions $f_\alpha$ are increasing and $f_\alpha(n)>n$ for
$n<\omega$. Fix a list $\langle \bar{m}_\xi:\xi < 2^{\aleph_0}\rangle$ of
all sequences $\bar{m} = \langle m_i:i < \omega\rangle$ such that $0 = m_0$
and $m_i +1 < m_{i+1}$.

For $\alpha<\gb$ we define:
\begin{enumerate}
\item[$(*)_1$]   $n_{\alpha,0} = 0$, $n_{\alpha,i+1} =
f_\alpha(n_{\alpha,i})$ (for $i<\omega$)  and $\bar{n}_\alpha= \langle
n_{\alpha,i}:i <\omega\rangle$;
\item[$(*)_2$]   $\bar{n}^0_\alpha =\langle 0,n_{\alpha,2},n_{\alpha,4},
  \ldots \rangle=\langle n^0_{\alpha,i}:i < \omega\rangle$ and
  $\bar{n}_\alpha^1 = \langle 0,n_{\alpha,3},n_{\alpha,5}, n_{\alpha,7},
  \ldots \rangle =\langle n^0_{\alpha,i}:i < \omega\rangle$.
\end{enumerate}
Observe that
\begin{enumerate}
\item[$(*)_3$] if $\bar{m}\in {}^\omega \omega$ is increasing, then for
  every large enough $\alpha<\gb$ we have:
\begin{enumerate}
\item[$(\alpha)$]  $(\exists^\infty i<\omega)(m_{i+1} < f_\alpha(m_i))$, and
  hence 
\item[$(\beta)$]  for at least one $\ell\in \{0,1\}$ we have 
\[\big(\exists^\infty i<\omega\big)\big(\exists j<\omega\big)
\big([m_i,m_{i+1})\subseteq [n^\ell_{\alpha,j}, n^\ell_{\alpha,j+1})\big).\]  
\end{enumerate}
\end{enumerate}
Now, for $\xi<2^{\aleph_0}$ we put:
\begin{enumerate}
\item[$(*)_4$] $\gamma(\xi)=\min\{\alpha<\gb: (\exists^\infty
  i<\omega)(f_\alpha(m_{\xi,i}) > m_{\xi,i+1})\}$;  
\item[$(*)_5$] $\ell(\xi) = \min\{\ell\le 1:(\exists^\infty
  i{<}\omega)(\exists j{<}\omega)([m_{\xi,i}, m_{\xi,i+1}) \subseteq
  [n^\ell_{\gamma(\xi),j}, n^\ell_{\gamma(\xi), j+1}))\}$;
\item[$(*)_6$] $\cU^1_\xi=\{i<\omega:(\exists j<\omega)
  ([m_{\xi,i},m_{\xi,i+1}) \subseteq  [n^{\ell(\xi)}_{\gamma(\xi),j},
  n^{\ell(\xi)}_{\gamma(\xi),j+1}))\}$.
\end{enumerate}
Note that $\gamma(\xi)$ is well defined by $(\alpha)$ of $(*)_3$, and so
also $\ell(\xi)$ is well defined (by $(\beta)$ of $(*)_3$). Plainly,
$\cU^1_\xi$ is an infinite subset of $\omega$. Now, for each
$\xi<2^{\aleph_0}$, we may choose $\cU^2_\xi$ so that 
\begin{enumerate}
\item[$(*)_7$]  $\cU^2_\xi\subseteq \cU^1_\xi$ is infinite and for any
  $i_1<i_2$ from $\cU^2_\xi$ we have  
\[(\exists j<\omega)(m_{\xi,i_1+1} <  n^{\ell(\xi)}_{\gamma(\xi),j}\ \&\ 
n^{\ell(\xi)}_{\gamma(\xi),j+1} < m_{\xi,i_2}).\]
\end{enumerate}
Let a function $g_\xi:\cU^2_\xi\longrightarrow \omega$ be such that 
\begin{enumerate}
\item[$(*)_8$] ${i\in \cU^2_\xi\ \&\ g_\xi(i)= j} \quad \Rightarrow \quad 
  {[m_{\xi,i},m_{\xi,i+1})\subseteq [n^{\ell(\xi)}_{\gamma(\xi),j},
  n^{\ell(\xi)}_{\gamma(\xi),j+1})}$.
\end{enumerate}
Clearly, $g_\xi$ is well defined and one-to-one.  (This is very important,
since it makes sure that the set $g_\xi[\cU^2_\xi]$ is infinite.) 
\medskip

Fix a sequence $\bar{B}=\langle B_{\zeta,t}:\zeta < \theta,\ t \in I_\zeta\rangle$
given by Lemma \ref{lemm} (so $\theta\leq\gb$ and $\bar{B}$ satisfies the
demands in (a)--(c) of \ref{lemm}).  By clause \ref{lemm}(c), for every
$\xi< 2^{\aleph_0}$, the set 
\[\big\{(\zeta,t):\zeta<\theta\mbox{ and }t \in I_\zeta\mbox{ and
}B_{\zeta,t} \cap g_\xi[\cU^2_\xi]\mbox{ is infinite }\big\}\]
has cardinality continuum. 
\medskip

Now, for each $\beta<\gb^+$ and $\xi<2^{\aleph_0}$ we choose a pair
$(\zeta_{\beta,\xi}, t_{\beta,\xi})$ such that 
\begin{enumerate}
\item[$(*)_9$]  $\zeta_{\beta,\xi}<\theta$ and $t_{\beta,\xi}\in
  I_{\zeta_{\beta,\xi}}$, 
\item[$(*)_{10}$]  $B_{\zeta_{\beta,\xi},t_{\beta,\xi}} \cap
  g_\xi[\cU^2_\xi]$ is infinite, and   
\item[$(*)_{11}$]  $t_{\beta,\xi}\notin \{t_{\alpha,\vare}:\vare<\xi$
or $\vare=\xi\ \&\ \alpha < \beta\}$.
\end{enumerate}
To carry out the choice we proceed by induction {\em first\/} on $\xi <
2^{\aleph_0}$, then on $\beta<\gb^+$. As there are $2^{\aleph_0}$ pairs
$(\zeta,t)$ satisfying clauses $(*)_9+(*)_{10}$ whereas clause $(*)_{11}$
excludes $\le (\gb^+ + |\xi|)\times\theta< 2^{\aleph_0}$ pairs (recalling
$\gb^+ < \frakg \le 2^{\aleph_0}$), there is such a pair at each stage
$(\beta,\xi)\in \gb^+\times 2^{\aleph_0}$.  

Lastly, for $\beta<\gb^+$ and $\xi<2^{\aleph_0}$ we let
\begin{enumerate}
\item[$(*)_{12}$] $\cU_{\beta,\xi}=g^{-1}_\xi[B_{\zeta_{\beta,\xi},
    t_{\beta,\xi}}] \cap \cU^2_\xi$
\end{enumerate}
(it is an infinite subset of $\omega$) and we put 
\begin{enumerate}
\item[$(*)_{13}$]  $A^+_{\beta,\xi} = \bigcup\{[m_{\xi,i},m_{\xi,i+1}):i\in
  \cU_{\beta,\xi}\}$, and  
\item[$(*)_{14}$]  $\cA_\beta = \{A \in [\omega]^{\aleph_0}$: for some
  $\xi < 2^{\aleph_0}$ we have $A \subseteq A^+_{\beta,\xi}\}$. 
\end{enumerate}

By the choice of $\langle \bar m_\xi:\xi<2^{\aleph_0}\rangle$,
$A^+_{\beta,\xi}$ and $\cA_\beta$ one easily verifies that for each
$\beta<\gb^+$: 
\begin{enumerate}
\item[$(*)_{15}$]  $\cA_\beta$ is a groupwise dense subset of
  $[\omega]^{\aleph_0}$. 
\end{enumerate}
Since we are assuming towards contradiction that $\frakg>\gb^+$, there is an
infinite $B \subseteq\omega$ such that 
\[(\forall \beta < \gb^+)(\exists A \in \cA_\beta)(B \subseteq^* A).\]
Hence for every $\beta<\gb^+$ we may choose $\xi(\beta)<2^{\aleph_0}$ such
that $B \subseteq^* A^+_{\beta,\xi(\beta)}$. Now, since $\gamma(\xi(\beta))
< \gb$ and $\zeta_{\beta,\xi(\beta)} < \theta \leq \gb$ and
$\ell(\xi(\beta)) \in \{0,1\}$, hence for some triple
$(\gamma^*,\zeta^*,\ell^*)$ we have that 
\begin{enumerate}
\item[$(\odot)_1$]  the set 
\[W =:\big\{\beta<\gb^+:\big(\gamma(\xi(\beta)), \zeta_{\beta,\xi(\beta)},  
\ell(\xi(\beta))\big)=\big(\gamma^*,\zeta^*,\ell^*\big) \big\}\]  
is unbounded in $\gb^+$.
\end{enumerate}
Note that if $\beta \in W$ then (recalling $(*)_{13}+(*)_8+(*)_{12}$)
\begin{enumerate}
\item[$(\odot)_2$]  $B \subseteq^* A^+_{\beta,\xi(\beta)}=\bigcup \big\{ 
  [m_{\xi(\beta),i}, m_{\xi(\beta),i+1}):i
  \in\cU_{\beta,\xi(\beta)}\big\}\subseteq$\\
$\bigcup\big\{[n^{\ell(\xi(\beta))}_{\gamma(\xi(\beta)),j}, n^{\ell(\xi(
  \beta))}_{\gamma(\xi(\beta)),j+1}):j=g_{\xi(\beta)}(i) \mbox{ for
  some }i \in \cU_{\beta,\xi(\beta)}\big\}\subseteq$\\ 
$\bigcup\big\{[n^{\ell(\xi(\beta))}_{\gamma(\xi(\beta)),j},
n^{\ell(\xi(\beta))}_{\gamma(\xi(\beta)),j+1}):j \in
B_{\zeta_{\beta,\xi(\beta)}, t_{\beta,\xi(\beta)}}\big\}$.
\end{enumerate}
Also, for $\beta \in W$ we have $\ell(\xi(\beta)) =\ell^*$,
$\gamma(\xi(\beta)) = \gamma^*$ and $\zeta(\beta,\xi(\beta))= \zeta^*$, so
it follows from $(\odot)_2$ that 
\begin{enumerate}
\item[$(\odot)_3$]  $B \subseteq^* \bigcup\big\{[n^{\ell^*}_{\gamma^*,j},
  n^{\ell^*}_{\gamma^*,j+1}):j \in B_{\zeta^*,t_{\beta,\xi(\beta)}}\}$ for
  every $\beta \in W$. 
\end{enumerate}

Consequently, if $\beta\ne \delta$ are from $W$, then the sets
\[\begin{array}{l}
\bigcup\big\{[n^{\ell^*}_{\gamma^*,j},n^{\ell^*}_{\gamma^*,j+1}):j \in 
B_{\zeta^*,t_{\beta,\xi(\beta)}}\big\}\mbox{ and}\\
\bigcup\big\{[n^{\ell^*}_{\gamma^*,j},n^{\ell^*}_{\gamma^*,j+1}):j \in
B_{\zeta^*,t_{\delta,\xi(\delta)}}\big\}
\end{array}\] 
are {\em not\/} almost disjoint. Hence, as $\langle
n^{\ell^*}_{\gamma^*,j}:j < \omega\rangle$ is increasing, necessarily the
sets $B_{\zeta^*,t_{\beta,\xi(\beta)}}$ and
$B_{\zeta^*,t_{\delta,\xi(\delta)}}$ are not almost disjoint.  So applying 
\ref{lemm}(b) we conclude that $t_{\beta,\xi(\beta)} =
t_{\delta,\xi(\delta)}$.  But this contradicts $\beta\ne \delta$ by
$(*)_{11}$, and we are done.
\end{proof}

\begin{definition}
\label{gf.1}
We define a cardinal characteristic $\frakg_\gf$ as the minimal cardinal
$\theta$ for which there is a sequence $\langle \cI_\alpha: \alpha< \theta
\rangle$ of groupwise dense {\em ideals\/} of $\cP(\omega)$ (i.e.,
$\cI_\alpha \subseteq [\omega]^{\aleph_0}$ is groupwise dense and $\cI_\alpha
\cup [\omega]^{< \aleph_0}$ is an ideal of subsets of $\omega$) such that  
\[\big(\forall B \in [\omega]^{\aleph_0}\big)\big(\exists\alpha<\theta\big)
\big(\forall A \in \cA_\alpha\big)\big(B\not\subseteq^* A\big).\]
\end{definition}

\begin{observation}
\label{gf.7}  $2^{\aleph_0} \ge \frakg_\gf\ge \frakg$.
\end{observation}

\begin{theorem}
\label{gf.14} 
$\frakg_\gf \le \gb^+$.
\end{theorem}

\begin{proof}
We repeat the proof of Theorem \ref{thmbg}. However, for $\beta< \gb^+$ the
family $\cA_\beta \subseteq [\omega]^{\le\aleph_0}$ does not have to be an
ideal.  So let $\cI_\beta$ be an ideal on $\cP(\omega)$ generated by
$\cA_\beta$ (so also $\cI_\beta$ is the ideal generated by
$\{A^+_{\beta,\xi}: \xi < 2^{\aleph_0}\} \cup [\omega]^{<\aleph_0}$).
Lastly, let $\cI'_\beta = \cI_\beta\setminus [\omega]^{< \aleph_0}$.

Assume towards contradiction that $B \in [\omega]^{\aleph_0}$ is such that
$(\forall \alpha<\gb^+)(\exists A \in\cI_\alpha)(B \subseteq^*A)$.  So for
each $\beta<\gb^+$ we can find $k_\beta < \omega$ and $\xi(\beta,0) <
\xi(\beta,1) < \ldots < \xi(\beta,k_\beta)<2^{\aleph_0}$ such that $B
\subseteq^* \bigcup\{A^+_{\beta,\xi(\beta,k)}:k \le k_\beta\}$.  Let $D$ be
a non-principal ultrafilter on $\omega$ to which $B$ belongs.  For each
$\beta<\gb^+$ there is $k(\beta)\le k_\beta$ such that
$A^+_{\beta,\xi(\beta,k(\beta))} \in D$.  As in the proof there for some
$(\gamma^*,\zeta^*,\ell^*,k^*,k(*))$ the following set is unbounded in
$\gb^+$: 
\[\begin{array}{rl}
W =: \big\{\beta < \gb^+:&k(\beta) = k(*),\ k_\beta = k^*, \  \gamma_{\xi(\beta,k(*))} =
  \gamma^*,\\
&\zeta_{\beta,\xi(\beta,k(*))} = \zeta^* \mbox{ and }
\ell(\xi(\beta,k(*)))=\ell^*\ \big\}.
\end{array}\]
As there it follows that:
\begin{enumerate}
\item[$(\odot)$]  if $\beta \in W$, then $\bigcup\big\{[
  n^{\ell^*}_{\gamma^*,j},n^{\ell^*}_{\gamma^*,j+1}):j \in
  B_{\zeta^*,t_{\beta,\xi(\beta,k(*))}}\big\}$ belongs to $D$.
\end{enumerate}
But for $\beta\ne \delta \in W$ those sets are almost disjoint whereas
$(\zeta^*,t_{\beta,\xi(\beta,k(*))}) \ne (\zeta^*,t_{\delta,\xi(\delta,k(*))})$
are distinct, giving us a contradiction.   
\end{proof}


\begin{thebibliography}{1}

\bibitem{BaSi89}
Bohuslav Balcar and Petr Simon.
\newblock {Disjoint refinement}.
\newblock In {\em Handbook of Boolean Algebras}, volume~2, pages 333--388.
  North-Holland, 1989.
\newblock Monk D., Bonnet R. eds.

\bibitem{BaJu95}
Tomek Bartoszy\'nski and Haim Judah.
\newblock {\em {Set Theory: On the Structure of the Real Line}}.
\newblock A K Peters, Wellesley, Massachusetts, 1995.

\bibitem{Bs89}
Andreas Blass.
\newblock {Applications of superperfect forcing and its relatives}.
\newblock In {\em {Set theory and its applications (Toronto, ON, 1987)}},
  volume 1401 of {\em Lecture Notes in Math.}, pages 18--40. Springer, Berlin,
  1989.

\bibitem{BsLa89}
Andreas Blass and Claude Laflamme.
\newblock {Consistency results about filters and the number of inequivalent
  growth types}.
\newblock {\em J. Symbolic Logic}, 54:50--56, 1989.

\bibitem{BrLo03}
J\"org Brendle and Maria Losada.
\newblock {The cofinality of the infinite symmetric group and groupwise
  density}.
\newblock {\em J. Symbolic Logic}, 68:1354--1361, 2003.

\bibitem{MdSh:843}
Heike Mildenberger and Saharon Shelah.
\newblock {Increasing the groupwise density number by c.c.c. forcing}.
\newblock {\em Annals of Pure and Applied Logic}.
\newblock math.LO/0404147.

\bibitem{MdSh:731}
Heike Mildenberger and Saharon Shelah.
\newblock {The relative consistency of ${\mathfrak g}<{\rm cf}({\rm
  Sym}(\omega))$}.
\newblock {\em Journal of Symbolic Logic}, 67:297--314, 2002.
\newblock math.LO/0009077.

\bibitem{Th98a}
Simon Thomas.
\newblock {Groupwise density and the confinality of the infinite symmetric
  group}.
\newblock {\em Arch. Math. Logic}, 37:483--493, 1998.

\bibitem{vD}
Eric~K. van Douwen.
\newblock {The integers and topology}.
\newblock In K.~Kunen and J.~E. Vaughan, editors, {\em Handbook of
  Set-Theoretic Topology}, pages 111--167. Elsevier Science Publishers, 1984.

\end{thebibliography}
\end{document}